\documentclass[a4paper,11pt]{amsart}
\usepackage{fullpage}
\usepackage{fontenc}
\usepackage{amsmath,amsfonts,amssymb}
\usepackage[sort]{natbib}
\makeatletter
\newcommand{\thorn}{{\fontencoding{T1}\selectfont\th}}
\newcommand{\@indepsymbol}[2]{#1\setbox0=\hbox{$#1x$}\kern\wd0\hbox to 0pt{\hss$#1\mid$\hss}\lower.9\ht0\hbox to 0pt{\hss$#1\smile$\hss}\kern\wd0}
\newcommand{\@nindepsymbol}[2]{#1\setbox0=\hbox{$#1x$}\kern\wd0\hbox to 0pt{\mathchardef
	\nn=12854\hss$#1\nn$\kern1.4\wd0\hss}\hbox to
	0pt{\hss$#1\mid$\hss}\lower.9\ht0 \hbox to
	0pt{\hss$#1\smile$\hss}\kern\wd0}
\newcommand{\ind}[1][]{\mathop{\mathpalette\@indepsymbol{}^{\!\!\!\!\rlap{$\scriptstyle\textnormal{#1}$}\,\,\,\,}}}
\newcommand{\nind}[1][]{\mathop{\mathpalette\@nindepsymbol{}^{\!\!\!\rlap{$\scriptstyle\textnormal{#1}$}\,\,\,}}}
\newcommand{\@Ind}[1][]{\mathpalette\@indepsymbol{}^{\!\!\!\!\mbox{$\scriptstyle\textnormal{#1}$}}}
\newcommand{\Ind}[1][]{\@Ind[\ \,]}
\newcommand{\newind}[4]{
	\newcommand{#1}{{\!\@Ind[#4]}}
	\newcommand{#2}{\ind[#4]}
	\newcommand{#3}{\nind[#4]}
}
\newind{\thInd}{\thind}{\nthind}{\thorn}
\renewcommand{\thInd}{\text{$\@Ind[\thorn]$\;}}
\makeatother

\newtheorem{thm}{Theorem}[section]
\newtheorem{cor}[thm]{Corollary}

\newtheorem{lem}[thm]{Lemma}

\newtheorem{fact}[thm]{Fact}
\theoremstyle{definition}
\newtheorem{defn}[thm]{Definition}
\theoremstyle{remark}

\theoremstyle{remark}

\theoremstyle{remark}

\title{Expansions which introduce no new open sets}
\author{Gareth Boxall and Philipp Hieronymi}
\date{\today}
\thanks{Gareth Boxall was supported by EPSRC grant EP/F009712/1. Philipp Hieronymi was supported by the Deutscher Akademischer Austausch Dienst.}

\begin{document}
\maketitle

\begin{abstract}
We consider the question of when an expansion of a topological structure has the property that every open set definable in the expansion is definable in the original structure. This question is related to and inspired by recent work of Dolich, Miller and Steinhorn on the property of having o-minimal open core. We answer the question in a fairly general setting and provide conditions which in practice are often easy to check. We give a further characterisation in the special case of an expansion by a generic predicate.
\end{abstract}

\section{Introduction}

Let $M$ be a many-sorted first-order topological structure in the sense of \cite{a}. So for each sort $S$ of $M$ there is a distinguished formula $\varphi_S(x,\bar{y})$ such that $x$ is of sort $S$ and the family of definable sets obtained by varying the parameters $\bar{y}$ forms a basis for a topology on $M_S$ (where $M_S$ is the set of realisations of the sort $S$). When we speak of a set being definable in $M$ we mean that it is a subset of some product $M_{S_1}\times...\times M_{S_n}$ and that it is first-order definable over some parameters. Each product $M_{S_1}\times...\times M_{S_n}$ is equipped with the appropriate product topology and so it makes sense to ask, for any definable set, whether or not it is open. Throughout this paper $M^*$ is an expansion of $M$ in which there are no new sorts. We consider the question ``when is it the case that every open set definable in $M^*$ is definable in $M$?''.

This question is inspired by recent work of Dolich, Miller and Steinhorn in \cite{acc} and \cite{acc2} on the property of having o-minimal open core. This property makes sense in the special case where $M$ is a one-sorted expansion of a model of the theory of dense linear orderings without endpoints and $\varphi(x,\bar{y})$ is ``$y_1<x<y_2$''. The open core of such an $M$ is a relational structure which has the same underlying set as $M$ and a predicate for each open set which is definable in $M$. The property of having o-minimal open core is an interesting generalisation of o-minimality which is studied extensively in \cite{acc}. There is also earlier work of Miller and Speissegger in \cite{cp}.

Two important classes of structures which are not o-minimal but which have o-minimal open core receive special attention in \cite{acc}. They are dense pairs of o-minimal ordered groups, as studied by van den Dries in \cite{l}, and expansions of o-minimal structures by a generic predicate, in the sense of Chatzidakis and Pillay in \cite{za}. In both cases the structure arises as an expansion of an o-minimal structure and it is proved in \cite{acc} that every open set definable in the expansion is definable in the original structure. So the question which we are addressing is intimately connected with the issue of having o-minimal open core. Our answer has the benefit of working outside the ordered setting.

We assume throughout that $M^*$ is $\kappa$-saturated and strongly $\kappa$-homogeneous for some sufficiently large $\kappa$. We fix a set $C$ in $M$ such that $|C|<\kappa$. It is clear that if $N^*\prec M^*$, $N$ is the reduct of $N^*$ to the language of $M$ and $C$ is a set in $N$ then the condition that every open set definable over $C$ in $N^*$ is definable over $C$ in $N$ is equivalent to the condition that every open set definable over $C$ in $M^*$ is definable over $C$ in $M$. Bearing this in mind it should be clear how to apply to $N$ and $N^*$ the results which we shall state for $M$ and $M^*$.

We denote by $B^{S_1...S_n}_{\bar{b}}$ the subset of $M_{S_1}\times...\times M_{S_n}$ defined by the formula $\varphi_{S_1}(x_1,\bar{b}_1)\wedge...\wedge\varphi_{S_n}(x_n,\bar{b}_n)$, where $\bar{b}=\bar{b}_1...\bar{b}_n$. We make use of the following two assumptions, the second of which is borrowed from Definition 4.6 of \cite{aae}.\\

\noindent Assumption (I): for any sorts $S_1,...,S_n$ of $M$ and $\bar{a}\in U\subseteq M_{S_1}\times...\times M_{S_n}$ such that $U$ is open, $\{\bar{b}:\bar{a}\in B^{S_1...S_n}_{\bar{b}}\subseteq U\}$ has non-empty interior.\\

\noindent Assumption (II): for any sorts $S_1,...,S_n$ of $M$ and $\bar{a}\in M_{S_1}\times...\times M_{S_n}$, $\{\bar{b}: \bar{a}\in B^{S_1...S_n}_{\bar{b}}\}$ is open. \\

Under assumption (I) we provide necessary and sufficient conditions for it to be the case that every open set definable over $C$ in $M^*$ is definable over $C$ in $M$. This is Theorem \ref{two}. Using somewhat similar methods we show that, under assumptions (I) and (II), it is enough to consider only open sets $U$ such that $U$ is equal to the interior of its closure. This is Theorem \ref{three}.

In Section 3 we apply Theorems \ref{two} and \ref{three} to classes of expansions which resemble the dense pairs of o-minimal ordered groups studied in \cite{l}, \cite{aca} and \cite{acc}. We answer a question from \cite{ae} about lovely pairs of dense o-minimal structures and prove that every open set definable in $(\mathbb{R},2^{\mathbb{Z}},2^{\mathbb{Z}}3^{\mathbb{Z}})$ is definable in $(\mathbb{R},2^{\mathbb{Z}})$. Outside the ordered setting we prove that every open set definable in a lovely pair of henselian valued fields of characteristic zero (and arbitrary residue field characteristic) is definable in the underlying henselian valued field.

In Section 4 we assume $M^*=M^g$ is an expansion of $M$ by a generic predicate in the sense of \cite{za}. It is established in \cite{acc} that if $M$ is an o-minimal expansion of a model of the theory of dense linear orderings without endpoints and $\varphi(x,\bar{y})$ is ``$y_1<x<y_2$'' then every open set definable in $M^g$ is definable in $M$. This is generalised in \cite{acc2} to the case where the assumption of o-minimality is replaced by the weaker assumption that $M$ has o-minimal open core. We use Theorem \ref{two} to give a characterisation of when it is the case that every open set definable over $C$ in $M^g$ is definable over $C$ in $M$, assuming $C=acl_M(C)$ and assumption (I). This is Theorem \ref{four}. We use it to reprove this result from \cite{acc2} and to obtain a similar result in which it is assumed that $M$ is a one-sorted geometric structure in the sense of Hrushovski and Pillay in \cite{ua}, open sets are infinite or empty and assumption (I) is true. We apply this work to henselian valued fields of characteristic zero considered first as one-sorted structures and then as two-sorted structures.

We would like to thank Anand Pillay and the audience at the mathematical logic seminar at Lyon 1 for some useful comments.

\section{General results}

For any sorts $S_1,...,S_n$ of $M$, let $A_{S_1...S_n}$ be the set of all $\bar{a}\in M_{S_1}\times...\times M_{S_n}$ such that the set of realisations of $tp_M(\bar{a}/C)$ is a dense subset of an open set. For any sorts $S_1,...,S_n$ of $M$, let $A^\prime_{S_1...S_n}$ be the set of all $\bar{a}\in A_{S_1...S_n}$ such that the set of realisations of $tp_{M^*}(\bar{a}/C)$ is dense in the set of realisations of $tp_M(\bar{a}/C)$. The following lemma is proved by a standard compactness argument.

\begin{lem}\label{one}
For any sorts $S_1,...,S_n$ of $M$ and any dense set $D\subseteq M_{S_1}\times...\times M_{S_n}$, $A_{S_1...S_n}\cap D$ is dense in $M_{S_1}\times...\times M_{S_n}$.
\end{lem}

\proof Suppose $S_1,...,S_n$ and $D\subseteq M_{S_1}\times...\times M_{S_n}$ are such that $D$ is dense in $M_{S_1}\times...\times M_{S_n}$ but $A_{S_1...S_n}\cap D$ is not dense in $M_{S_1}\times...\times M_{S_n}$. Let $U\subseteq M_{S_1}\times ...\times M_{S_n}$ be a non-empty open set such that $A_{S_1...S_n}\cap D\cap U=\emptyset$. Enumerate as $(p_\alpha)_{\alpha<\lambda}$ the complete types over $C$ in the sense of $M$ which are realised in $D\cap U$. It follows that, for each $\alpha<\lambda$, the set of realisations of $p_\alpha$ is nowhere dense. We obtain a decreasing sequence $(B^{S_1...S_n}_{\bar{b}_\alpha})_{\alpha<\lambda}$ of non-empty basic open subsets of $U$ as follows. Let $\alpha<\lambda$ and suppose $B^{S_1...S_n}_{\bar{b}_\beta}$ has already been chosen for all $\beta<\alpha$. By compactness there is some $B^{S_1...S_n}_{\bar{b}_{\alpha}^\prime}\subseteq \bigcap_{\beta<\alpha}B^{S_1...S_n}_{\bar{b}_\beta}\cap U$ such that $B^{S_1...S_n}_{\bar{b}_{\alpha}^\prime}\neq \emptyset$. Given that the set of realisations of $p_{\alpha}$ is not dense in $B^{S_1...S_n}_{\bar{b}_{\alpha}^\prime}$, choose $B^{S_1...S_n}_{\bar{b}_\alpha}\subseteq B^{S_1...S_n}_{\bar{b}_{\alpha}^\prime}$ such that $B^{S_1...S_n}_{\bar{b}_{\alpha}}\neq \emptyset$ and $p_\alpha$ is not realised in $B^{S_1...S_n}_{\bar{b}_\alpha}$. By compactness, $\bigcap_{\alpha<\lambda}B^{S_1...S_n}_{\bar{b}_\alpha}$ has non-empty interior. Therefore $\emptyset\neq D\cap \bigcap_{\alpha<\lambda}B^{S_1...S_n}_{\bar{b}_\alpha}\subseteq D\cap U$. However, for all $\alpha<\lambda$, $p_\alpha$ is not realised in $D\cap \bigcap_{\alpha<\lambda}B^{S_1...S_n}_{\bar{b}_\alpha}$. This is a contradiction. \endproof

\begin{thm}\label{two}
Suppose assumption (I) is true. The following are equivalent:

(1) for any sorts $S_1,...,S_n$ of $M$, $A^\prime_{S_1...S_n}$ is dense in $M_{S_1}\times...\times M_{S_n}$,

(2) for any sorts $S_1,...,S_n$ of $M$, $A^\prime_{S_1...S_n}=A_{S_1...S_n}$,

(3) every open set definable over $C$ in $M^*$ is definable over $C$ in $M$.
\end{thm}

\proof Suppose (1). Let $U$ be an open set definable over $C$ in $M^*$. Then $U\subseteq M_{S_1}\times...\times M_{S_n}$ for some sorts $S_1,...,S_n$ of $M$. Let $\bar{a}\in U$. There is some $B^{S_1...S_n}_{\bar{b}}$ such that $\bar{a}\in B^{S_1...S_n}_{\bar{b}}\subseteq U$. Let $S^\prime_1,...,S^\prime_m$ be such that $\bar{b}\in M_{S^\prime_1}\times...\times M_{S^\prime_m}$. By assumption (I) we may assume $\bar{b}$ is an interior point of $\{\bar{b}^\prime:\bar{a}\in B^{S_1...S_n}_{\bar{b}^\prime}\subseteq U\}$. By (1) we may further assume $\bar{b}\in A^\prime_{S^\prime_1...S^\prime_m}$. Let $\sigma$ be an automorphism of the structure $M$ which fixes $C$ pointwise. Then $\sigma(\bar{b})$ is an interior point of $\{\bar{b}^\prime:\sigma(\bar{a})\in B^{S_1...S_n}_{\bar{b}^\prime}\}$. Since $\bar{b}\in A^\prime_{S^\prime_1...S^\prime_m}$ it follows that there is some $\bar{b}^\prime\models tp_{M^*}(\bar{b}/C)$ such that $\sigma(\bar{a})\in B^{S_1...S_n}_{\bar{b}^\prime}$. Since $\bar{b}^\prime\models tp_{M^*}(\bar{b}/C)$, $B^{S_1...S_n}_{\bar{b}^\prime}\subseteq U$. Therefore $U$ is invariant under automorphisms of the structure $M$ which fix $C$ pointwise. Since $U$ is definable in $M^*$ and $M^*$ is sufficiently saturated, it follows that $U$ is definable over $C$ in $M$. Therefore (3).

Suppose not (2). Let $S_1,...,S_n$ be such that $A^\prime_{S_1...S_n}\neq A_{S_1...S_n}$. Let $\bar{a}\in A_{S_1...S_n}\setminus A^\prime_{S_1...S_n}$. Let $W$ be an open set such that the set of realisations of $tp_M(\bar{a}/C)$ is a dense subset of $W$. Let $B^{S_1...S_n}_{\bar{b}}$ be a non-empty basic open subset of $W$ in which $tp_{M^*}(\bar{a}/C)$ is not realised. By compactness there is a set $X$ definable over $C$ in $M^*$ such that $\bar{a}\in X$ and $X\cap B^{S_1...S_n}_{\bar{b}}=\emptyset$. Let $U$ be the interior of $M_{S_1}\times...\times M_{S_n}\setminus X$. Then $U$ is definable over $C$ in $M^*$. However $U$ is not definable over $C$ in $M$, since $tp_M(\bar{a}/C)$ is realised in $B^{S_1...S_n}_{\bar{b}}\subseteq U$ while $\bar{a}\notin U$. Therefore not (3).

(2)$\Rightarrow$(1) by Lemma \ref{one}. \endproof

The following is something of a definable set (as opposed to type) version of Theorem \ref{two}.

\begin{thm}\label{three}
Suppose assumptions (I) and (II) are true. Suppose every open set which is definable over $C$ in $M^*$ and equal to the interior of its closure is definable over $C$ in $M$. Then every open set which is definable over $C$ in $M^*$ is definable over $C$ in $M$.
\end{thm}

\proof Let $U\subseteq M_{S_1}\times...\times M_{S_n}$ be an open set definable over $C$ in $M^*$. Let $Y=\{\bar{b}:B^{S_1...S_n}_{\bar{b}}\subseteq U\}$. Let $V$ be the interior of $Y$. It follows from assumption (I) that $U$ is definable over $C$ in $M$ provided $V$ is definable over $C$ in $M$. It remains only to show that $V$ is equal to the interior of its closure. For this it is sufficient to show that $\overline{V}\subseteq Y$ (where $\overline{V}$ is the closure of $V$). Suppose $\overline{V}\nsubseteq Y$. Let $\bar{b}\in \overline{V}\setminus Y$. Then $B^{S_1...S_n}_{\bar{b}}\nsubseteq U$. Let $\bar{a}\in B^{S_1...S_n}_{\bar{b}}\setminus U$. By assumption (II) there is some $\bar{b}^\prime\in V\subseteq Y$ such that $\bar{a}\in B^{S_1...S_n}_{\bar{b}^\prime}\subseteq U$. This is a contradiction.  \endproof

It is clear that the saturation and homogeneity assumptions are not required for Theorem \ref{three}.

\section{Examples which resemble dense pairs}

We begin this section with two corollaries. Theorem 5.2 of \cite{ap} is in a similar spirit to these results, but it assumes that $M$ is a dense o-minimal ordered group. The following is a corollary of Theorem \ref{two}.

\begin{cor}\label{cor1}
Suppose assumption (I) is true. Suppose that for any sorts $S_1,...,S_n$ of $M$ there is a set $D_{S_1...S_n}\subseteq M_{S_1}\times...\times M_{S_n}$ such that the following conditions are satisfied:

(1) $D_{S_1...S_n}$ is dense in $M_{S_1}\times...\times M_{S_n}$,

(2) for every $\bar{a}\in D_{S_1...S_n}$ and every open set $U\subseteq M_{S_1}\times ...\times M_{S_n}$, if $tp_M(\bar{a}/C)$ is realised in $U$ then $tp_M(\bar{a}/C)$ is realised in $U\cap D_{S_1...S_n}$,

(3) for every $\bar{a}\in D_{S_1...S_n}$, $tp_{M^*}(\bar{a}/C)$ is implied by $tp_M(\bar{a}/C)$ in conjunction with ``$\bar{a}\in D_{S_1...S_n}$''.

Then every open set definable over $C$ in $M^*$ is definable over $C$ in $M$.
\end{cor}

\proof For any sorts $S_1,...,S_n$ of $M$, let $\bar{a}\in A_{S_1...S_n}\cap D_{S_1...S_n}$. Let $W$ be an open set such that the set of realisations of $tp_M(\bar{a}/C)$ is a dense subset of $W$. Let $U\subseteq W$ be a non-empty open subset. Then $tp_M(\bar{a}/C)$ is realised in $U$. By (2) there is some $\bar{a}^\prime\in U\cap D_{S_1...S_n}$ such that $\bar{a}^\prime\models tp_M(\bar{a}/C)$. By (3), $\bar{a}^\prime\models tp_{M^*}(\bar{a}/C)$. Therefore $\bar{a}\in A^\prime_{S_1...S_n}$. So $A_{S_1...S_n}\cap D_{S_1...S_n}\subseteq A^\prime_{S_1...S_n}$. It follows by (1) and Lemma \ref{one} that $A^\prime_{S_1...S_n}$ is dense in $M_{S_1}\times...\times M_{S_n}$. The result follows by Theorem \ref{two}. \endproof

The following is a corollary of Theorem \ref{three}.

\begin{cor}\label{cor}
Suppose assumptions (I) and (II) are true. Suppose that for any sorts $S_1,...,S_n$ of $M$ there is a set $D_{S_1...S_n}\subseteq M_{S_1}\times...\times M_{S_n}$ such that the following conditions are satisfied:

(1) $D_{S_1...S_n}$ is dense in $M_{S_1}\times...\times M_{S_n}$,

(2) for every $X\subseteq M_{S_1}\times...\times M_{S_n}$ which is definable over $C$ in $M$, either $X$ is nowhere dense or $X$ has non-empty interior,

(3) for every $X\subseteq M_{S_1}\times...\times M_{S_n}$ which is definable over $C$ in $M^*$, there exists $Y\subseteq M_{S_1}\times ...\times M_{S_n}$ which is definable over $C$ in $M$ and such that $X\cap D_{S_1...S_n}=Y\cap D_{S_1...S_n}$.

Then every open set definable over $C$ in $M^*$ is definable over $C$ in $M$.
\end{cor}

\proof For any sorts $S_1,...,S_n$ of $M$, let $U\subseteq M_{S_1}\times...\times M_{S_n}$ be open and definable over $C$ in $M^*$. By (3) there is a set $Y\subseteq M_{S_1}\times...\times M_{S_n}$ which is definable over $C$ in $M$ and such that $U\cap A_{S_1...S_n}\cap D_{S_1...S_n}=Y\cap A_{S_1...S_n}\cap D_{S_1...S_n}$. Let $V$ be the interior of $\overline{Y}$. Let $Z$ be either $V\setminus Y$ or $Y\setminus V$. Suppose $Z\cap A_{S_1...S_n}\cap D_{S_1...S_n}\neq\emptyset$. Let $\bar{a}\in Z\cap A_{S_1...S_n}\cap D_{S_1...S_n}$. Since $Z$ is definable over $C$ in $M$, it follows that the set of all realisations of $tp_M(\bar{a}/C)$ is contained in $Z$. Therefore $Z$ is not nowhere dense. Since $Z$ has empty interior, this contradicts (2). Therefore $V\cap A_{S_1...S_n}\cap D_{S_1...S_n}=Y\cap A_{S_1...S_n}\cap D_{S_1...S_n}=U\cap A_{S_1...S_n}\cap D_{S_1...S_n}$. By (1) and Lemma \ref{one}, $A_{S_1...S_n}\cap D_{S_1...S_n}$ is dense in $M_{S_1}\times...\times M_{S_n}$. Since $U$ and $V$ are both open, it follows that $\overline{U}=\overline{V}$. By Theorem \ref{three} we may assume that $U$ is the interior of $\overline{U}$. Then $U$ is definable over $C$ in $M$. \endproof

For the rest of this section we assume that $M$ and $M^*$ are one-sorted. We adapt our notation accordingly. For example, $A_{S_1...S_n}$ becomes $A_n$. Algebraic closure in the sense of $M$ is denoted by $acl_M$. We recall the following definition from \cite{ua}.

\begin{defn}\label{geom}
$M$ is geometric if

(1) $M$ is infinite,

(2) $Th(M)$ eliminates $\exists^\infty$ and

(3) $acl_M$ is a pregeometry on $M$.
\end{defn}

The second condition is sometimes called ``uniform finiteness''. Suppose $M$ is geometric and $M^*$ is a lovely pair of models of $Th(M)$. So $M^*$ expands $M$ by a single unary predicate $P$ such that $P(M)\prec M$ as structures in the language of $M$ and, for every finite $A\subseteq M$ and every infinite $X\subseteq M$ such that $X$ is definable in $M$, $X\cap P(M)\neq\emptyset$ and $X\nsubseteq acl_M(P(M)\cup A)$. This situation (without the topological aspect which we are assuming here) is investigated by Berenstein and Vassiliev in \cite{ae}. Among other things, it is known to generalise the dense pairs of o-minimal ordered groups studied in \cite{l}. The following corollary generalises a result in \cite{acc}.

\begin{cor}\label{pair}
Suppose $M$ is geometric and the expansion $M^*$ is a lovely pair of models of $Th(M)$. Suppose assumption (I) is true and every open set is infinite or empty. Suppose $C$ is $acl_M$-independent from $P(M)$ over $C\cap P(M)$. Then every open set definable over $C$ in $M^*$ is definable over $C$ in $M$.
\end{cor}

\proof For all $n<\omega$, let $D_n=\{\bar{a}\in M^n: \bar{a}$ is $acl_M$-independent over $P(M)\cup C\}$. Since non-empty basic open sets are infinite and definable in $M$, it follows from the definition of lovely pairs and compactness that $D_n$ is dense in $M^n$. This gives condition (1) of Corollary \ref{cor1}. Let $\bar{a}\in D_n$ and let $U$ be an open set such that $tp_M(\bar{a}/C)$ is realised in $U$. Let $\bar{b}\in M^m$ and $\bar{a}^\prime\in M^n$ be such that $\bar{a}^\prime\in B_{\bar{b}}^n\subseteq U$ and $\bar{a}^\prime \models tp_M(\bar{a}/C)$. By assumption (I), compactness and the fact that non-empty open sets are infinite, we may assume $\bar{b}$ is $acl_M$-independent over $C\bar{a}^\prime$. It follows, by condition (3) of Definition \ref{geom}, that $\bar{a}^\prime$ is $acl_M$-independent over $C\bar{b}$. Then, by the definition of lovely pairs and compactness, $tp_M(\bar{a}^\prime/C\bar{b})$ is realised in $D_n$. This gives condition (2) of Corollary \ref{cor1}. Condition (3) of Corollary \ref{cor1} is proved in \cite{ae}. \endproof

The case of Corollary \ref{pair} in which $M$ expands a model of the theory of dense linear orderings without endpoints and $\varphi(x,\bar{y})$ is ``$y_1<x<y_2$'' answers a question from \cite{ae}. It is proved in \cite{acc} that if $M^*$ is a dense pair of o-minimal ordered groups in the sense of \cite{l} then every open set definable in $M^*$ is definable in $M$. Suppose $M$ is a henselian valued field of characteristic zero (by which we mean that the field has characteristic zero and the residue field has arbitrary characteristic). Suppose the language of $M$ is the language of rings together with a unary predicate for the valuation ring $\mathcal{O}$. As in \cite{a}, take $\varphi(x,\bar{y})$ to be ``$\frac{x-y_1}{y_2}\in \mathcal{O}$''. It is known, and follows from work of van den Dries in \cite{l2}, that $M$ is geometric.  It is clear that assumption (I) is true and open sets are infinite or empty. So Corollary \ref{pair} applies and we conclude that every open set definable in a lovely pair of models of $Th(M)$ is definable in the underlying henselian valued field.

Some other examples of expansions which resemble dense pairs are referred to in \cite{acc}. These include $(\mathbb{R},\mathbb{U})$, the real field expanded by a predicate for the group of complex roots of unity (studied by Zilber in \cite{b}), and $(\mathbb{R},2^{\mathbb{Z}}3^{\mathbb{Z}})$, the real field expanded by a predicate for the group $2^{\mathbb{Z}}3^{\mathbb{Z}}$ (studied by van den Dries and G\"{u}nayd\i n in \cite{la}). These examples are considered in greater generality in \cite{ob} and \cite{la}, but we do not wish to go into the details of that here. It is known that any open set definable in $(\mathbb{R},\mathbb{U})$ or $(\mathbb{R},2^{\mathbb{Z}}3^{\mathbb{Z}})$ is definable in $\mathbb{R}$ (see \cite{aca}, \cite{acc} and \cite{ap}). We reprove this here. Let $N=\mathbb{R}$ and $N^*=(N,P(N))$ where $P(N)=\mathbb{U}$ or $P(N)=2^{\mathbb{Z}}3^{\mathbb{Z}}$. Suppose $M^*$ elementarily extends $N^*$ and $M$ is the reduct of $M^*$ to the language of $N$. Suppose $C\prec N^*$ and $C$ is countable. For all $n<\omega$, let $D_n=\{\bar{a}\in M^n: \bar{a}$ is $acl_M$-independent over $P(M)\cup C\}$. We check the conditions of Corollary \ref{cor}. Condition (1) is true for $N$ in place of $M$ and $D_n$ restricted to $N^n$, by the well known fact that the complement of any countable subset of $\mathbb{R}$ is dense in $\mathbb{R}$. It follows, using compactness, that condition (1) is true for $M$ and $D_n$. Condition (2) is a well known property of real-closed fields. Condition (3) is proved in \cite{b} and \cite{la}. So the conclusion of Corollary \ref{cor} follows and it remains true when $M$ is replaced by $N$ and $M^*$ is replaced by $N^*$.

A more complicated example involves taking $N=(\mathbb{R},2^{\mathbb{Z}})$ and $N^*=(\mathbb{R},2^{\mathbb{Z}},2^{\mathbb{Z}}3^{\mathbb{Z}})$. This expansion is studied by G\"{u}nayd\i n in \cite{ayh}, the structure $(\mathbb{R},2^{\mathbb{Z}})$ having previously been studied by van den Dries in \cite{l3}. Let $P(N)=2^{\mathbb{Z}}3^{\mathbb{Z}}$. Suppose $M^*$ elementarily extends $N^*$ and $M$ is the reduct of $M^*$ to the language of $N$. Let $M^\prime$ be the reduct of $M$ to the language of $\mathbb{R}$. Suppose $C\prec N^*$ and $C$ is countable. For all $n<\omega$, let $D_n=\{\bar{a}\in M^n: \bar{a}$ is $acl_{M^\prime}$-independent over $P(M)\cup C\}$. We check the conditions of Corollary \ref{cor}. Condition (1) is true by the same reasoning as in the previous example. Condition (2) follows from the quantifier elimination proved in \cite{l3} (see \cite{chris}). Condition (3) is proved in \cite{ayh}. The conclusion of Corollary \ref{cor} follows and so every open set definable in $(\mathbb{R},2^{\mathbb{Z}},2^{\mathbb{Z}}3^{\mathbb{Z}})$ is definable in $(\mathbb{R},2^{\mathbb{Z}})$. This answers a question asked by Chris Miller and discussed at a meeting at the Fields Institute in Toronto in January 2009. A positive answer was expected at the time. The example was of interest as an instance of the phenomenon of a non-d-minimal structure which has d-minimal but not o-minimal open core. We refer the reader to \cite{chris} for background and more on d-minimality.

\section{Adding a generic predicate}

We drop the assumption that $M$ is one-sorted. Throughout this section we assume that $M^*=M^g$ is an expansion of $M$ by a generic predicate as defined by Chatzidakis and Pillay in \cite{za}. Let $G(M)$ be the set of realisations of this generic predicate. So there is an infinite set $X$ such that $X$ is $\emptyset$-definable in $M$, $G(M)\subseteq X$ and the following two facts, which we recall from \cite{za}, are true.

\begin{fact}\label{za1} Let $n<\omega$ and let $a_1...a_n=\bar{a}\in X^n$ be such that, for all $i\in \{1,...,n\}$, $a_i\notin acl_M(C)$. Let $J\subseteq \{1,...,n\}$. Then there is some $\bar{a}^\prime=a_1^\prime...a_n^\prime$ such that $\bar{a}^\prime\models tp_M(\bar{a}/C)$ and $a^\prime_i\in G(M)$ if and only if $i\in J$.
\end{fact}

\begin{fact}\label{za2}
Let $\bar{a}$ and $\bar{a}^\prime$ be tuples from $M$ of length less than $\kappa$. Suppose there is an automorphism $\sigma$ of the structure $M$ such that $\sigma$ fixes $C$ pointwise, $\sigma(\bar{a})=\bar{a}^\prime$ and, for all $b\in X\cap acl_M(C\bar{a})$, $b\in G(M)$ if and only if $\sigma(b)\in G(M)$. Then $tp_{M^*}(\bar{a}/C)=tp_{M^*}(\bar{a}^\prime/C)$.
\end{fact}

Recall from \cite{alf} that a consistent formula $\psi(\bar{y},\bar{d})$ in the language of $M$ strongly divides over $C$ if $\bar{d}\nsubseteq acl_M(C)$ and $\{\psi(\bar{y},\bar{d}^\prime):\bar{d}^\prime\models tp_M(\bar{d}/C)\}$ is $k$-inconsistent for some $k<\omega$. Using compactness it is observed in \cite{alf} that, in this case, there is a set $Z$ which is definable over $C$ in $M$ and such that $\bar{d}\in Z$ and $\{\psi(\bar{y},\bar{d}^\prime):\bar{d}^\prime\in Z\}$ is $k$-inconsistent.

\begin{thm}\label{four}
Suppose $C=acl_M(C)$ and assumption (I) is true. The following are equivalent:

(1) for every $n<\omega$, $\bar{d}\in X^n$ and formula $\psi(\bar{y},\bar{d})$ in the language of $M$, if $\psi(\bar{y},\bar{d})$ defines a non-empty open set then $\psi(\bar{y},\bar{d})$ does not strongly divide over $C$,

(2) every open set definable over $C$ in $M^g$ is definable over $C$ in $M$.
\end{thm}

\proof Suppose not (2). By Theorem \ref{two} there are sorts $S_1,...,S_n$ of $M$ such that $A^\prime_{S_1...S_n}\neq A_{S_1...S_n}$. Let $\bar{a}\in A_{S_1...S_n}\setminus A^\prime_{S_1...S_n}$. Let $W$ be an open set such that the set of realisations of $tp_M(\bar{a}/C)$ is a dense subset of $W$. Let $B^{S_1...S_n}_{\bar{b}}\subseteq W$ be non-empty and such that $tp_{M^*}(\bar{a}/C)$ is not realised in $B^{S_1...S_n}_{\bar{b}}$. Let $(\hat{a}_i)_{i<\lambda}$ be an enumeration of $X\cap acl_M(C\bar{a})\setminus C$. By Fact \ref{za2}, for every $\bar{a}^\prime\hat{a}_0^\prime\hat{a}_1^\prime...\models tp_M(\bar{a}\hat{a}_0\hat{a}_1.../C)$ such that $\bar{a}^\prime\in B^{S_1...S_n}_{\bar{b}}$, there is some $i<\lambda$ such that $\hat{a}_i\in G(M)$ if and only if $\hat{a}^\prime_i\notin G(M)$. It follows by Fact \ref{za1} and compactness that, for every $\bar{a}^\prime\hat{a}_0^\prime\hat{a}_1^\prime...\models tp_M(\bar{a}\hat{a}_0\hat{a}_1.../C)$ such that $\bar{a}^\prime\in B^{S_1...S_n}_{\bar{b}}$, there is some $i<\lambda$ such that $\hat{a}^\prime_i\in acl_M(C\bar{b})$. Therefore, by compactness, there exist a subtuple $(i_1,...,i_m)$ of $\lambda$ and a tuple $(D_1,...,D_m)$ of finite subsets of $X$ such that, for every $\bar{a}^\prime\hat{a}_0^\prime\hat{a}_1^\prime...\models tp_M(\bar{a}\hat{a}_0\hat{a}_1.../C)$ such that $\bar{a}^\prime\in B^{S_1...S_n}_{\bar{b}}$, there is some $j\in \{1,...,m\}$ such that $\hat{a}^\prime_{i_j}\in D_j$.

Keeping $(i_1,...,i_m)$ fixed, we may assume $|D_1|$ is minimal and that, subject to this, $|D_2|$ is minimal and so on. We may assume $\bar{b}$ has been chosen so that this remains true even when $\bar{b}$ is replaced by $\bar{b}^\prime$ for some $\emptyset\neq B^{S_1...S_n}_{\bar{b}^\prime}\subseteq B^{S_1...S_n}_{\bar{b}}$. We may also assume each $D_j$ is non-empty. It is easy to check that $(D_1,...,D_m)$ is uniquely determined by $(i_1,...,i_m)$, $(|D_1|,...,|D_m|)$ and the fact that, for every $\bar{a}^\prime\hat{a}_0^\prime\hat{a}_1^\prime...\models tp_M(\bar{a}\hat{a}_0\hat{a}_1.../C)$ such that $\bar{a}^\prime\in B^{S_1...S_n}_{\bar{b}}$, there is some $j\in \{1,...,m\}$ such that $\hat{a}^\prime_{i_j}\in D_j$. By compactness there is a formula $\psi^\prime(\bar{y},\bar{z})$ over $C$ in the language of $M$ such that, for every $\bar{b}^\prime$ such that $\emptyset\neq B^{S_1...S_n}_{\bar{b}^\prime}\subseteq B^{S_1...S_n}_{\bar{b}}$, $\psi^\prime(\bar{b}^\prime,\bar{z})$ defines the set of all enumerations of $D_1\cup...\cup D_m$. Let $\bar{d}$ be such an enumeration. By assumption (I), the set defined by $\psi^\prime(\bar{y},\bar{d})$ has non-empty interior. Let $\psi(\bar{y},\bar{z})$ be a formula over $C$ in the language of $M$ such that $M\models\psi(\bar{y},\bar{z})\rightarrow \psi^\prime(\bar{y},\bar{z})$ and $\psi(\bar{y},\bar{d})$ defines the interior of the set defined by $\psi^\prime(\bar{y},\bar{d})$. Let $k$ be the number of permutations of $\bar{d}$. We may assume that, for every tuple $\bar{b}^\prime$, $\psi(\bar{b}^\prime,\bar{z})$ has at most $k$ realisations. Since $\hat{a}_{i_j}\notin C=acl_M(C)$ for every $j\in\{1,...,m\}$, we have $\bar{d}\nsubseteq acl_M(C)$ and clearly $\{\psi(\bar{y},\bar{d}^\prime):\bar{d}^\prime\models tp_M(\bar{d}/C)\}$ is $k$-inconsistent. Therefore not (1).

Suppose not (1). Suppose $d_1...d_n=\bar{d}\in X^n$ and $\psi(\bar{y},\bar{d})$ defines a non-empty open set and strongly divides over $C$. We may assume $d_1\notin acl_M(C)$. Let $Z\subseteq X^n$ be definable over $C$ in $M$ and such that $\bar{d}\in Z$, $\psi(\bar{y},\bar{d}^\prime)$ defines a non-empty open set for every $\bar{d}^\prime\in Z$ and $\{\psi(\bar{y},\bar{d}^\prime):\bar{d}^\prime\in Z\}$ is $k$-inconsistent for some $k<\omega$. We may assume $Z$ and $k$ are such that $k$ is minimal. Using Fact \ref{za1} it follows that $\bigvee\{\psi(\bar{y},\bar{d}^1)\wedge...\wedge\psi(\bar{y},\bar{d}^{k-1}):\bar{d}^1,...,\bar{d}^{k-1}$ are distinct members of $Z$ and $d^1_1,...,d^{k-1}_1\in G(M)\}$ defines a non-empty union of non-empty open sets which is therefore non-empty and open. Let $U$ be this union. So $U$ is definable over $C$ in $M^g$. Let $\bar{d}^1,...,\bar{d}^{k-1}$ be distinct realisations of $tp_M(\bar{d}/C)$ such that $\{\psi(\bar{y},\bar{d}^1),...,\psi(\bar{y},\bar{d}^{k-1})\}$ is consistent and $d^1_1,...,d^{k-1}_1\in G(M)$. By Fact \ref{za1} there is some $\bar{d}^{\prime 1}...\bar{d}^{\prime k-1}\models tp_M(\bar{d}^1...\bar{d}^{k-1}/C)$ such that $d^{\prime 1}_1\notin G(M)$. It follows that the non-empty set defined by $\psi(\bar{y},\bar{d}^1)\wedge...\wedge\psi(\bar{y},\bar{d}^{k-1})$ is contained in $U$ and the non-empty set defined by $\psi(\bar{y},\bar{d}^{\prime 1})\wedge...\wedge\psi(\bar{y},\bar{d}^{\prime k-1})$ is disjoint from $U$. So $U$ is not definable over $C$ in $M$. Therefore not (2). \endproof

The following is a result of Dolich, Miller and Steinhorn from \cite{acc2}.

\begin{cor}\label{dms}
Suppose $M$ is a one-sorted expansion of a model of the theory of dense linear orderings without endpoints and $\varphi(x,\bar{y})$ is ``$y_1<x<y_2$''. If $M$ has o-minimal open core then every open set definable in $M^g$ is definable in $M$.
\end{cor}

\proof Assume $C=acl_M(C)$. Suppose condition (1) of Theorem \ref{four} is not satisfied. Let $\bar{d}\in X^n$ and let $\psi(\bar{y},\bar{d})$ be a formula in the language of $M$ such that $\psi(\bar{y},\bar{d})$ defines a non-empty open set and strongly divides over $C$. Let $Z\subseteq X^n$ be definable over $C$ in $M$ and such that $\bar{d}\in Z$, $\psi(\bar{y},\bar{d}^\prime)$ defines a non-empty open set for every $\bar{d}^\prime \in Z$ and $\{\psi(\bar{y},\bar{d}^\prime):\bar{d}^\prime \in Z\}$ is $k$-inconsistent for some $k<\omega$. We may assume $Z$ and $k$ are such that $k$ is minimal. Then the set defined by $\bigvee\{\psi(\bar{y},\bar{d}^1)\wedge...\wedge\psi(\bar{y},\bar{d}^{k-1}):\bar{d}^1,...,\bar{d}^{k-1}$ are distinct members of $Z\}$ is definable over $C$ in $M$ and is an infinite union of pair-wise disjoint non-empty open sets each of which is definable in $M$. This contradicts the fact that $M$ has o-minimal open core. \endproof

The following corollary generalises the o-minimal special case of the previous one.

\begin{cor}\label{rank1}
Suppose $M$ is geometric (as in Definition \ref{geom}, so $M$ is one-sorted). Suppose assumption (I) is true and every open set is infinite or empty. Then every open set definable in $M^g$ is definable in $M$.
\end{cor}

\proof Assume $C=acl_M(C)$. Suppose condition (1) of Theorem \ref{four} is not satisfied. Let $\bar{d}\in X^n$ and let $\psi(\bar{y},\bar{d})$ be a formula in the language of $M$ such that $\psi(\bar{y},\bar{d})$ defines a non-empty open set and strongly divides over $C$. By compactness and the fact that non-empty open sets are infinite, there is a tuple $\bar{b}$ such that $M\models \psi(\bar{b},\bar{d})$ and $\bar{b}$ is $acl_M$-independent over $C\bar{d}$. Since $\psi(\bar{y},\bar{d})$ strongly divides over $C$, $\bar{d}\subseteq acl_M(C\bar{b})$ by Remark 3.2 of \cite{clifalf}. Since $\bar{d}\nsubseteq acl_M(C)$, this contradicts condition (3) of Definition \ref{geom}. \endproof

If $M$ is a henselian valued field of characteristic zero (and arbitrary residue field characteristic), considered as a one-sorted structure as in Section 3, then it follows by Corollary \ref{rank1} that every open set definable in $M^g$ is definable in $M$. Now suppose $M$ is the same henselian valued field but this time considered as a two-sorted structure: one sort $S_1$ for the field, equipped with the language of rings, one sort $S_2$ for the value group together with $\infty$, equipped with the language of ordered groups, and the valuation map $v$ between the two sorts. The topology on $M_{S_1}$ is as before only this time we take the formula $\varphi_{S_1}(x,\bar{y})$ to be ``$v(x-y_1)>y_2$'', where $x$ and $y_1$ are of sort $S_1$ and $y_2$ is of sort $S_2$. We take $\varphi_{S_2}(x,\bar{y})$ to be ``$y_1<x<y_2$'', where all three variables are of sort $S_2$. It is clear that assumption (I) is true for this structure.

Let $\psi(\bar{y},\bar{d})$ be a formula in the language of $M$ such that $\bar{d}\in M_{S_1}^n$ for some $n<\omega$, $\psi(\bar{y},\bar{d})$ defines a non-empty open set and $\psi(\bar{y},\bar{d})$ strongly divides over $C=acl_M(C)$. We may assume $\bar{y}=y_1...y_ky_{k+1}...y_m$ where $y_1,...,y_k$ are of sort $S_1$ and $y_{k+1},...,y_m$ are of sort $S_2$. Let $\psi^\prime(\bar{y}^\prime,\bar{d})$ be a formula which defines the set of all $b_1...b_kb_{k+1}^\prime...b_m^\prime$ such that $M\models \psi(b_1...b_kv(b_{k+1}^\prime)...v(b_m^\prime),\bar{d})$. Since $v$ is continuous, $\psi^\prime(\bar{y}^\prime,\bar{d})$ defines an open set. Let $C_1$ and $C_2$ be such that $C_1\subseteq M_{S_1}$, $C_2\subseteq M_{S_2}$ and $C=C_1\cup C_2$. Let $C_2^\prime\subseteq M_{S_1}$ be such that $|C_2^\prime|<\kappa$ and $v(C_2^\prime)=C_2$. Let $C^\prime=C_1\cup C_2^\prime$. It is known that the two-sorted structure which we are considering has the one-sorted structure, which we were considering, as a reduct on the sort $S_1$ and that the two-sorted structure is $\emptyset$-interpretable in the one-sorted structure. So $\psi^\prime(\bar{y}^\prime,\bar{d})$ is equivalent to a formula in the one-sorted structure. Suppose $\psi(\bar{y},\bar{d})$ strongly divides over $C$ in the two-sorted structure. It is clear that then $\psi^\prime(\bar{y}^\prime,\bar{d})$ strongly divides over $C$ in the two-sorted structure. We may assume $C^\prime_2$ has been chosen so that $\bar{d}$ is not algebraic over $C^\prime$ in either structure. It follows that $\psi^\prime(\bar{y}^\prime,\bar{d})$ strongly divides in the one-sorted structure over the algebraic closure of $C^\prime$ in the one-sorted structure. This contradicts Theorem \ref{four} for the one-sorted structure. Therefore condition (1) of Theorem \ref{four} is true for the two-sorted structure, provided $X$ is a subset of a power $M_{S_1}$. So no new open sets are introduced when expanding the field sort by a generic predicate.


\bibliographystyle{abbrv}

\end{document}